\magnification=\magstep1 
\overfullrule=0pt
\def\eqde{\,{\buildrel \rm def \over =}\,} 
 \def\la{\lambda}  \def\ga{\gamma} 
\def\J{{\cal J}}    \def\c{{\cal C}}   \def\i{{\rm i}} \def\B{{\cal B}}
\def\si{\sigma} \def\eps{\epsilon}    
\def\s{{\cal S}}       \def\I{{\cal I}}  
  \def\L{{\Lambda}} \def\M{{\cal M}}
\def\E{{\cal E}}  \def\p{{\cal P}}  \def\D{{\cal D}} 
 
\font\huge=cmr10 scaled \magstep2
\def\QED{\vrule height6pt width6pt depth0pt}
\font\smcap=cmcsc10     
  
\def\ade{${\cal ADE}_7$-type invariant}   
\input amssym.def
\def\Z{{\Bbb Z}}  \def\Q{{\Bbb Q}}

\def\boxit#1{\vbox{\hrule\hbox{\vrule{#1}\vrule}\hrule}}
\def\splus{\,\,{\boxit{$+$}}\,\,}
\def\stimes{\,\,{\boxit{$\times$}}\,\,}

{\nopagenumbers
\rightline{August, 1998}
\bigskip\bigskip
\centerline{{\bf \huge  The Level 2 and 3 Modular Invariants}}\bigskip
\centerline{{\bf \huge  for the Orthogonal Algebras}}
\bigskip \bigskip   \centerline{Terry Gannon}\medskip
\centerline{{\it Department of Mathematical Sciences, University of Alberta,}}
\centerline{{\it Edmonton, Alberta, Canada, T6G 1G8} } \smallskip
\centerline{{e-mail: tgannon@math.ualberta.ca}}
\bigskip\medskip

\noindent{\bf Abstract.}
This paper finds for each affine algebra $B_r^{(1)}$ and $D_r^{(1)}$
all modular invariant 1-loop partition functions at level $\le 3$.
Previously, only those at level 1 were classified. 
An extraordinary number of exceptionals appear at level 2
-- indeed this is the motivation for this paper -- and we find infinitely
many new ones there. The only level 3 exceptionals occur for
$B_2^{(1)}\cong C_2^{(1)}$ and $D_7^{(1)}$, and the latter appear to be new.
The $B_{2,3}$ and $D_{7,3}$ exceptionals are cousins of the
${\cal E}_6$-exceptional and ${\cal E}_8$-exceptional, respectively, of
$A_1^{(1)}$, while the level 2 exceptionals are related to the lattice
invariants of affine $u(1)$.


\vfill\eject}\pageno=1

\noindent{{\bf 1. Introduction}}\bigskip 

\noindent
Over the past decade or so, much work has been directed towards one aspect
of the classification of conformal field theories: the classification of
modular invariant partition functions corresponding to the affine algebras.
This classification question can be asked for any choice of algebra $X_r^{(1)}$
and level $k\in\Z_+\eqde\{0,1,2,\ldots\}$. An elaborate machinery has been
developed,
and we can be optimistic about our chances for the complete classification,
at least when $X_r$ is simple. Nevertheless, at the present time few of these
classifications have been accomplished: the main successes are $X_r=A_1$,
$A_2$, and $U_1\oplus\cdots\oplus U_1$ (the quasi-rational unitary case
should be interpreted rationally by equating representations with equal character),
for all levels $k$ [2,6,8]; and levels $k\le 3$ for all $A_r$ [7].

The result for $A_1$ falls into an A-D-E pattern [2]. Philippe
Ruelle [15] discovered a connection between Jacobians of Fermat curves,
and the $A_2$ classification. There is a natural relationship [8] between
rational points on Grassmannians and the $U_1\oplus\cdots\oplus U_1$ classification.
Relations between these partition functions and subfactors in von 
Neumann algebra theory
is discussed e.g.\ in [3]. For these reasons, as well of course for
conformal field theory itself, it should certainly be of interest to make
further efforts to obtain complete lists of these partition functions,
and to understand better the curious relationships between those lists and
other areas of mathematics and mathematical physics.

The finitely many characters $\chi_\la$, $\la\in P_+$, of the integrable
highest weight representations of $X_r^{(1)}$ at level $k$, carry a
representation of the modular group SL$_2(\Z)$ [12]. The matrices
$\left(\matrix{0&-1\cr 1&0\cr}\right)$ and $\left(\matrix{1&1\cr 0&1}\right)$,
which generate SL$_2(\Z)$, get sent to the Kac-Peterson matrices  $S$ and
$T$, respectively. The entries $S_{\mu\nu}$ are related to values of
Lie group characters at elements of finite order, while those of the
diagonal matrix $T$ are related to the eigenvalues of the quadradic Casimir. Our
classification problem reduces
to finding all matrices $M=(M_{\mu\nu})_{\mu,\nu\in P_+}$ which obey

\smallskip\item{(P1)} $MS=SM$ and $MT=TM$

\item{(P2)} $M_{\mu\nu}\in \Z_+$ for all $\mu,\nu\in P_+$

\item{(P3)} $M_{00}=1$\smallskip

\noindent (P1) says that
the corresponding partition function should be a modular invariant, (P2)
applies because the entries of $M_{\mu\nu}$ count certain `primary fields',
while (P3) says that the vacuum is unique. Any such matrix $M$ is called a
{\it physical invariant}.

An approach has slowly evolved to handle these classifications. It breaks the
problem into 2 parts: \smallskip

\item{(1)} find all possible values of $M_{\mu 0}$ and $M_{0 \mu}$,
for all $\mu\in P_+$;

\item{(2)} find all physical invariants $M$, for each possible choice of
values $M_{\mu 0}$ and $M_{0 \mu}$ found in (1). \smallskip

\noindent{The} point is that the values in (1) are severely constrained.
 What is found is that for any
physical invariant of almost every choice of $X_{r,k}$ (i.e.\ algebra $X_r^{(1)}$
at level $k$),
$$M_{\mu 0}\ne 0\quad {\rm or}\quad M_{0 \mu}\ne 0
\qquad \Longrightarrow\qquad\mu\in\s (0)\ ,\eqno(1.1)$$
where $\s$ is the group of symmetries of the Coxeter-Dynkin diagram of 
$X_r^{(1)}$ ($\s$ acts on $P_+$ by permuting the indices $0,1,2,\ldots,r$
of the weights $\la$). The orbit $\s (0)$ in (1.1) is the set of simple currents.
For example, for $A_{1,k}$, all but two
of the physical invariants obey (1.1) (those two are the so-called
$\E_6$ and $\E_8$ exceptionals, at $k=10$ and 28 respectively). Any
physical invariant obeying (1.1) is called an \ade, by analogy with the
$A_{1,k}$ classification. All evidence points to the validity of the
following conjecture: \medskip

\centerline{{\it for any fixed choice of  simple algebra
$X_r$,}}
\centerline{{\it all
but finitely many physical invariants of $X_{r,k}$ will obey (1.1).}}\medskip 

\noindent Hence an important -- and recently completed [9,10] -- component
of the classification of all physical invariants for $X_{r,k}$, is to
find all \ade s. This can be thought of as the generic situation in (2).
By contrast, little general work has been done on (1).

However this general programme breaks down for the orthogonal algebras at level
2. They behave completely uncharacteristically: part (1) has far too many
solutions and (1.1) rarely holds. Thus  $B_{r,2}$ and $D_{r,2}$ must be
treated using novel arguments. This is the primary motivation for this paper.

On the other hand, $B_{r,3}$ and $D_{r,3}$ behave generically (all their
physical invariants are of ${\cal ADE}_7$-type, except for $B_{2,3}$ and
$D_{7,3}$),
and we include their classification as an indication of more normal
behaviour. There are only three level 3 exceptionals. 

Level $k\le 3$ physical invariants are also known for the algebra $A_r^{(1)}$
[7]. Low-level classifications essentially reduce to low-rank ones,
because of a curious relation called {\it rank-level duality} (see e.g.\
[14]). For example, the Kac-Peterson matrices of $\widehat{so}_n$ level
$k$ and $\widehat{so}_k$ level $n$ are related. In particular, we find that
$B_{r,2}$ and $D_{r,2}$ are related to $U_1$ at levels $2r+1$ and $2r$,
respectively, while $B_{r,3}$ and $D_{r,3}$ are related to $A_{1,4r+4}$
and $A_{1,4r+2}$.

\bigskip\bigskip\noindent{{\bf 2. The list of physical invariants}} 

\bigskip \noindent We begin this
section with a quick  review, and then we list all physical invariants for
the orthogonal algebras at $k\le 3$.

\bigskip\noindent{{\it 2.1. Definitions.}}\quad
A highest-weight $\la=\la_0\L_0+\la_1\L_1+\cdots+\la_r\L_r\in P_+$ in
$B_{r,k}$ ($r\ge 3)$ satisfies $k=\la_0+\la_1+2\la_2+\cdots+2\la_{r-1}
+\la_r$, while for $D_{r,k}$ ($r\ge 4$) $\la\in P_+$ satisfies
$k=\la_0+\la_1+2\la_2+\cdots+2\la_{r-2}+\la_{r-1}+\la_r$; in both cases
the $\L_i$ are the fundamental weights and all $\la_i\in\Z_+$. For $B_{r,k}$
put $n=k+2r-1$, while for $D_{r,k}$ put $n=k+2r-2$. The Weyl vector
is $\rho=\sum \L_i$.  Usually we will drop the redundant $\L_0$ component
of the weights. Note that $B_{2,k}$ is more properly called $C_{2,k}$,
and so can, should, and will be ignored in the following.

The symmetries $\s$ of the Coxeter-Dynkin diagrams
will play a major role. Those fixing
the 0th node are called {\it conjugations}, while others -- called {\it
simple currents} -- form an abelian group we will call $\s_{sc}$.

$B_{r,k}$ has no nontrivial conjugation.
For any $D_{r,k}$, there is a conjugation $C_1$ interchanging $\la_{r-1}
\leftrightarrow\la_r$.
Put $C_0=I$, the identity. When $r=4$, there are four
additional conjugations $C_2,\ldots,C_5$ -- 
these six $C_i$ for $D_4^{(1)}$ correspond to the different
permutations of its Dynkin labels $\la_1,\la_3,\la_4$.

$B_{r,k}$ has a simple current $J_b$
of order 2, given by $J_b\la=(\la_1,\la_0,\la_2,\ldots,\la_r)$.
There are three non-trivial simple currents for $D_{r,k}$, namely $J_v$, $J_s$
and $J_c=J_v \circ J_s$, defined by
$J_v\la=(\la_1,\la_0,\la_2,\ldots,\la_{r-2},\la_r,\la_{r-1})$ and
$$J_s\la=\cases{
(\la_r,\la_{r-1},\la_{r-2},\ldots,\la_1,\la_0) & if $r$ is even \cr
(\la_{r-1},\la_r,\la_{r-2},\ldots,\la_1,\la_0) & if $r$ is odd\cr}\ .$$ 
Write $\J_b=\{id.,J_b\}$, $\J_v=\{id.,J_v\}$, $\J_s=\{id.,J_s\}$, and
$\J_d=\{id.,J_v,J_s,J_c\}$.
By a {\it spinor} for $B_{r,k}$ or $D_{r,k}$, respectively, is
meant any weight $\la\in P_+$ with $\la_r$ or $\la_{r-1}+\la_r$ odd.
Write $\p_b$ and $\p_v$ for the sets of nonspinors.

We say $\la$ and $\mu$ are $M$-{\it coupled} if either $M_{\la \mu}\ne 0$
or $M_{\mu \la}\ne 0$. By a {\it  positive invariant} we mean a matrix $M$ 
commuting with the
corresponding Kac-Peterson matrices $S$ and $T$, with the additional
property that each $M_{\mu \nu}\ge 0$. By a {\it physical invariant}, we
mean a positive invariant with each $M_{\mu \nu}\in\Z$, and obeying (P3).
By an {\it \ade}, we mean a
physical invariant $M$ satisfying (1.1). Finally, by an {\it automorphism
invariant}, we mean a physical invariant obeying
$$M_{\la 0}=M_{0 \la}=\delta_{\la,0}\ .\eqno(2.1a)$$
Automorphism invariants are important examples of physical invariants. It turns
 out (see Lemma 3.1(c)) that any automorphism invariant will be a permutation
matrix, i.e.\ there will be a permutation $\pi$ of $P_+$ such that
$$M_{\la \mu}=\delta_{\mu,\pi\la}\ .\eqno(2.1b)$$

Each conjugation $C$ defines an automorphism invariant,
which we will also denote by $C$,
obtained by taking $\pi=C$ in (2.1b). Moreover, the matrix products $C
\, M$ and $M\,C$ of $C$ with any other physical invariant $M$ will also be a
physical invariant.

The primary reason for the importance of simple currents  is: let $J\in
\s_{sc}$, then 
$$S_{J\mu,\nu}=\exp[2\pi \i\,Q_J(\nu)]\,S_{\mu\nu}\ ,\eqno(2.2a)$$
for some number $Q_J(\nu)$ [13,17]. $Q_{J_b}(\mu)=\la_r/2$, while
$Q_{J_v}(\mu)=(\la_{r-1}+\la_r)/2$ and
$$ Q_s(\la)=\sum_{j=1}^{r-2}j\la_j/2-{r-2 \over 4}\la_{r-1}
-{r \over 4}\la_r\ .$$
The matrix $T$ also behaves similarly under $\s_{sc}$:
$${(J\mu+\rho)^2 -(\mu+\rho)^2\over 2n}\equiv 
{R(J)\,(N-1)\over 2N}-Q_J(\mu)\quad({\rm mod}\ 1)\ ,\eqno(2.2b)$$
where $N$ is the order of $J$, and where $R(J)$ is some integer.
$R(J_b)=R(J_v)=2k$ and $R(J_s)=R(J_c)=N_s\,(N_s-1)kr/4$, where                     
$N_s$ is the order of $J_s$: $N_s=2$ or 4 depending on whether or not $r$ is
even.

 From these equations, it is possible to find a sequence of physical 
(in fact ${\cal ADE}_7$-type) invariants, for
 each $J\in\s_{sc}$. In particular, define [17]
$$\I[J]_{\mu,\nu}=\sum_{\ell=1}^{N}
\delta_{J^{\ell}\mu,\nu}\,\delta^1\bigl(Q_J(\mu)+{\ell\over 2N}
R(J)\bigr)\ ,\eqno(2.3)$$
where $\delta^1(x)=1$ if $x\in\Z$ and $=0$ otherwise. 
For example, $\I[id.]=I$, the identity matrix. $\I[J]$ will be
a physical invariant iff $R(J)$ is even.

Any physical invariant not constructable in these
standard ways out of simple currents and conjugations  is called an
{\it exceptional invariant}.

\medskip\noindent{{\it 2.2. The list of physical invariants for $B_{r,2}$.}}
Here $n=2r+1$. $P_+$ consists of precisely
$r+4$ weights: $0$, $J_b0=2\L_1$, $\L_r$, $J_b\L_r=\L_1+\L_r$, $\ga^i\eqde
\L_i$ for $i<r$,        and $\ga^r\eqde 2\L_r$.  Write
$\ga^0$ for the weight 0.
To minimise subscripts, we will usually abbreviate `$J_b$' to `$J$'.
Extra exceptionals exist when  $n$ is a perfect square: $\sqrt{n}\in\Z$.
In that case $4|r$, and it is convenient to introduce the following notation:
when $8|r$, write $\la^r\eqde \L_r$ and $\mu^r\eqde J\L_r$; otherwise write
$\la^r\eqde J\L_r$ and $\mu^r\eqde \L_r$. Also write ${\cal C}=\{\ga^a\ne
0\,|\,\sqrt{n}\ {\rm divides}\ a\}$.

Define matrices $\B(d,\ell)$, $\B(d_1,\ell_1|d_2,\ell_2)$, $\B^i$, $\B^{ii}$,
$\B^{iii}$, $\B^{iv}$ by:
$$\eqalignno{\B(d,\ell)_{J^i\ga^a,J^i\ga^b}=&\left\{\matrix{2&{\rm if}\
d|a,\  d|b, {\rm and\ both}\ a\ne 0,\ b\ne 0\cr
0&{\rm if}\ n\not|da\ {\rm or}\ b\not\equiv\pm a\ell\ ({\rm mod}\ 
d)\cr 1& {\rm otherwise}\cr}\right. &\cr
\B(d,\ell)_{J^i\L_r,J^i\L_r}=&\,1&\cr}$$
and all other entries are 0, where $a,b\in\{0,1,\ldots,r\}$ and $i\in\{0,1\}$;
$$\eqalignno{\B(d_1,\ell_1|d_2,\ell_2)=&\,
{1\over 2}(\B(d_1,\ell_1)+\B(d_2,\ell_2))\,\I[J_b]&\cr
\B^{i}_{00}=\B^{i}_{0\ga}=\B^{i}_{\ga 0}=&\,\B^i_{\ga\ga'}=\B^i_{\la^r\ga}=
\B^{i}_{\ga\la^r}=\B^{i}_{\mu^r\mu^r}=\B^{i}_{\la^r,J0}=
\B^{i}_{J0,\la^r}=1&\cr
\B^{ii}_{00}=\B^{ii}_{0\ga}=\B^{ii}_{\ga 0}=&\,\B^{ii}_{\ga\ga'}=
\B^{ii}_{0\la^r}=\B^{ii}_{\la^r 0}=
\B^{ii}_{\la^r\la^r}=\B^{ii}_{\la^r\ga}=\B^{ii}_{\ga\la^r}=1&\cr}$$
and all other entries are 0, where $\ga,\ga'\in{\cal C}$. Finally,
$\B^{iii}=\B^i\,\I[J_b]$ and $\B^{iv}=\I[J_b]\,\B^i$.

In section 4.1 we will prove:

\medskip\noindent{\smcap Theorem 2.1}.\quad {\it Let $M$ be a physical invariant of
$B_{r,2}$. Then $M$ equals one of the following:\smallskip

\item{(a)} $\B(d,\ell)$ for any divisor $d$ of $n=2r+1$ obeying $n|d^2$,
and for any integer $0\le \ell<{d^2\over 2n}$ obeying $\ell^2\equiv 1$
(mod ${d^2\over n}$);\smallskip

\item{(b)} $\B(d_1,\ell_1|d_2,\ell_2)$ for any divisors $d_i$ of $n$ obeying
$n|d_i^2$, and for any integers $0\le \ell_i< {d_i^2\over 2n}$ obeying
$\ell_i^2\equiv 1$ (mod ${d_i^2\over n}$);\smallskip

\item{(c)} when $n$ is a perfect square, there are 4 remaining
physical invariants: $\B^i$, $\B^{ii}$, $\B^{iii}${,} and $\B^{iv}$.}\medskip

The only redundancy here is that $\B(d_1,\ell_1|d_2,\ell_2)=\B(d_2,
\ell_2|d_1,\ell_1)$.    There are a total of $D$ distinct $\B(d,\ell)$'s, where $D$
is the number of divisors $d'\le \sqrt{n}$ of $n$. All but one of these,
namely $\B(n,1)=I$, are exceptional.
There are precisely $D\,(D+1)/2$ distinct physical invariants in part (b),
and all but one of them (namely $\B(n,1|n,1)=\I[J_b]$) are exceptional. When
$\sqrt{n}\in\Z$, the physical invariants $\B^i$, $\B^{ii}$, $\B^{iii}$,
$\B^{iv}$ are all distinct and exceptional. Of all $B_{r,2}$ physical
invariants, only $\B^{iii}$ and $\B^{iv}$ are not symmetric matrices.

Most of these exceptionals are new. The $\B(d,\ell)$ for the special case
$d=n$ are the exceptional automorphism invariants found in [4]. The
$\B(d,\ell)$ for the special case $\ell=1$ are the exceptional invariants
given in [16]. Multiplying these gives all invariants of type $\B(d,\ell)$.

For example, when $3\le r\le 12$, respectively, there are precisely
2, 9, 2, 2, 5, 2, 2, 5, 2, 9 physical invariants for $B_{r,2}$. All nine
$B_{4,2}$ physical invariants can be found in the Appendix B of [18], and the
correspondence between his notation and ours is: 
$Z_1=\B(9,1)$, $Z_2=\B(9,1|9,1)$, $Z_3=\B(3,1)$, $Z_4=\B^{i}$, $Z_5=\B^{ii}$,
$Z_6=\B(3,1|3,1)$, $Z_7=\B^{iii}$, $Z_8=\B^{iv}$, and $Z_9=\B(3,1|9,1)$. Verstegen
attributed this $B_{4,2}$ richness to the existence of conformal embeddings
involving $E_8, A_8$ and $D_8$, and so was unaware that $B_{4,2}$ is merely the
tip of an iceberg!

\medskip\noindent{{\it 2.3. The list of physical invariants for $D_{r,2}$.}}
There are $r+7$ weights: 0, $2\L_1=J_v0$, $2\L_{r-1}=J_c0$, $2\L_r=J_s0$,
$\L_r$, $\L_1+\L_{r-1}=J_v\L_r$, $\L_{r-1}$, $\L_1+\L_r=J_v\L_{r-1}$,
$\la^i\eqde \L_i$ for $1\le i\le r-2$, and $\la^{r-1}\eqde
\L_{r-1}+\L_r$. Write $n=2r$. Write
$\la^0$ for the weight $0$ and $\la^r$ for $2\L_r$, and $J$ for $J_v$.
Additional exceptionals occur when $r$ is an even perfect square, and 
in this case write ${\cal C}_j=\{\la^b\ne 0\,|\,2{b\over \sqrt{r}}\equiv \pm j\
 ({\rm mod}\ 8)\}$ for $j=0,1,2,3,4$. 

Define the matrices $\D(d,\ell)$, $\D(d_1,\ell_1|d_2,\ell_2)$, $\D^i$, $\D^{ii}$,
and $\D^{iii}$, as follows:
$$\eqalignno{\D(d,\ell)_{J^i\la^a,J^i\la^b}=&\left\{\matrix{2&{\rm if}\
d|a,\ d|b,\ 2d|(a+b),\ {\rm and}\ \{a,b\}\subseteq\{1,\ldots,r-1\}\cr
0&{\rm  if}\ r\not| da\ {\rm or}\ 
b\not\equiv \pm a\ell\ ({\rm mod}\ 2d)\cr 1&{\rm otherwise}}\right. &\cr
\D(d,\ell)_{\la_s\la_s}=&\left\{\matrix{1&{\rm if}\ 2d\not|
r\cr 2&{\rm if}\ \la_s\in\{\L_r,\L_1+\L_{r-1}\}\ {\rm and}\ 2d|r\cr
0&{\rm otherwise}}\right.&\cr}$$
and all other entries are 0,
where $a,b\in\{0,1,\ldots,r\}$, $i\in\{0,1\}$, and $\la_s$ is any spinor;
$$\eqalignno{\D(d_1,\ell_1|d_2,\ell_2)=&\,
{1\over 2}(\D(d_1,\ell_1)+\D(d_2,\ell_2))\,\I[J_v]&\cr
\D^{i}_{J^j\L_r,J^j\L_r}=&\,\D^{i}_{\L_r\mu}=\D^{i}_{\mu\L_r}=
\D^{i}_{J\L_r,\mu'}=\D^{i}_{\mu',J\L_r}=\D^{i}_{\la\la'}
=\D^{i}_{\ga\ga'}&\cr=&\,\D^{i}_{J'0,J''0}=\D^{i}_{J'0,\nu}=
\D^{i}_{\nu,J'0}=\D^{i}_{J'J0,\nu'}=\D^{i}_{\nu',J'J0}=1&\cr}$$
 where $\la,\la'\in{\cal C}_0\cup{\cal C}_4$, $\mu\in 
{\cal C}_1$, $\mu'\in{\cal C}_3$, $\ga,\ga'\in{\cal C}_2$, $\nu\in
{\cal C}_0$, $\nu'\in{\cal C}_4$, $J',J''\in\J_s$, and $j\in\{0,1\}$.
All other entries  equal 0.  Finally, $\D^{ii}=\D^i\,\I[J_v]$ and $\D^{iii}
=\I[J_v]\,\D^{i}$.

In section 4.2 we prove:

\medskip\noindent{\smcap Theorem 2.2}.\quad {\it Let $M$ be a physical invariant of
$D_{r,2}$. Then $M$ equals one of the following, for arbitrary conjugations
$C_i,C_j$:

\item{(a)} $C_i\,\D(d,\ell)\,C_j$ for  any divisor $d$ of $r$
 obeying $r|d^2$,
and for any integer $1\le \ell\le{d^2\over r}$ obeying $\ell^2\equiv 1$
(mod ${4d^2\over r}$); \smallskip

\item{(b)} $\D(d_1,\ell_1|d_2,\ell_2)$ for any divisors $d_i$ of $r$ obeying
$r|d_i^2$, as well as the additional property that $2d_1|r$ iff $2d_2|r$,
and for any integers $1\le \ell_i\le {d_i^2\over r}$ obeying
$\ell_i^2\equiv 1$ (mod ${4d_i^2\over r}$);

\smallskip\item{(c)} when $r$ is a perfect square and $16|r$, there are 8 
other physical invariants: $C_i\,\D^i\,C_j${,} $C_i\,\D^{ii}$, and
$\D^{iii}\,C_j$.       }\medskip

Take $C_i=I$ in (a) unless $2d|r$. To count these physical invariants,
define $D\{m\}$ to be the number of divisors $d\le\sqrt{m}$ of $m$, and
let $2^c$ be the exact power of 2 dividing $r$. Put $D=D\{r\}$ and $D_1=D-D_0$,
where $D_0=D\{2r/2^c\}$ when 4 divides $r$, and $D_0=D$ otherwise. 
For $r\ne 4$, there are precisely $2D_0+4D_1$ distinct physical invariants
in (a) and $D_0(D_0+1)/2+D_1(D_1+1)/2$ in (b).
For $D_{4,2}$, there are a total of 16 physical invariants (namely, 6 of the
form $\D(4,1)\, C_i$,
9 of the form $C_i\, \D(2,1)\, C_j$, and $\D(2,1|2,1)$).

Most of these physical invariants are new, and all but the conjugations of
$\D(r,1)=I$ and $\D(r,1|r,1)=\I[J_v]$ (for all $r$), and $\D(r,r-1)=\I[J_s]$
and $\D(r,r-1|r,r-1)=\I[J_v]\,\I[J_s]$ (when $r\equiv 2$ (mod 4)), and
$\D({r\over 2},1)=\I[J_s]$ and $\D({r\over 2},1|{r\over 2},1)=\I[J_v]\,
\I[J_s]$ (when $4|r$) are exceptional. The
exceptionals $\D(d,\ell)$ for the special cases $d=r$ and $\ell=1$, respectively,
first appeared in [4] and [16].

For example, for $4\le r\le 16$, respectively, there are precisely 16, 3, 8, 3, 8, 7,
7, 3, 12, 3, 7, 7, and 22 physical invariants. Of these, 0, 0, 2, 0, 0, 4, 1,
0, 4, 0, 1, 4, and 14 are exceptional.

\medskip\noindent{{\it 2.4. The orthogonal algebras at levels 3 and 1.}}
Next consider $B_{r,3}$, $r\ge 3$. In section 5.1 we will show that there are
no exceptional physical invariants (incidently, there is one for $C_{2,3}$
-- recall $C_2\cong B_2$):

\medskip\noindent{\smcap Theorem 2.3}.\quad {\it The only physical invariants for
$B_{r,3}$ are the identity matrix $M=I$ and the simple current invariant
$M=\I[J_b]$.}\medskip

The physical invariant $\I[J_b]$ was first found in [1]. It is an
automorphism invariant; the associated permutation is order 2. The $C_{2,3}$
exceptional is a conformal embedding.

Next look at $D_{r,3}$, for $r\ge 4$. In section 5.2 we will show that the
only exceptional invariants occur at $D_{7,3}$.

\medskip\noindent{\smcap Theorem 2.4}.\quad {\it The complete list of $D_{r,3}$ physical
 invariants  is:

 \item{$\bullet$} \quad $C_i$ and $C_i\,\I[J_v]$, valid for all $k$;

\item{$\bullet$} \quad in addition for $r\equiv 4$ (mod $8$), $C_i\,M$ for
$M=\I[J_s]$, $\I[J_c]$, $\I[J_s]\,\I[J_c]$ and $\I[J_c]\,\I[J_s]$;

\item{$\bullet$} \quad in addition when $8$ divides $r$, $C_i\,M\,C_j$ for
$M=\I[J_s]$ and $\I[J_v]\,\I[J_s]$;

\item{$\bullet$}\quad finally, for $D_{7,3}$, there are the exceptionals
$C_i\,\E(D_{7,3})$, where
$${\cal E}(D_{7,3})\eqde
\sum_{J\in\J_d}|\chi_{J0}+\chi_{J(\L_1+\L_5)}|^2+\sum_{J\in\J_d}|
\chi_{J\L_3}+\chi_{J(\L_1+\L_6+\L_7)}|^2\ ,$$
\noindent where $C_i$ and $C_j$ are any conjugations (i.e.\ $i,j\in\{0,1\}$
for $r\ne 4$, and $i,j\in \{0,1,\ldots,5\}$ for $r=4$).}\medskip

$\I[J_v]$ and $\I[J_s]$ were first found in [1]. To our knowledge,
$\E(D_{7,3})$ has never appeared before in the literature. $\I[J_v]$ is
an automorphism invariant;
the permutation will be order 2. When $r\equiv 4$ (mod 8), $\I[J_s]$ and
$\I[J_c]$ will
also be order 2 automorphism invariants. When 8 divides $r$, $\I[J_s]$
will be a direct sum of $r+2$ matrices $\left(\matrix{1&1\cr 1&1}\right)$
and $2r+4$ matrices $(0)$.

For completeness, let us repeat the level 1 classification [5].
There are no exceptionals. For $B_{r,1}$ there is only the identity matrix
$I$. For $D_{r,1}$, when 4 does not divide $r$, there are only 2 invariants:
$I$ and $C_1=\I[J_v]$. When 4 divides $r$, there are a total of 6 invariants:
$C_i$ and $C_i\,\I[J_s]\,C_j$, for $i,j\in\{0,1\}$.

\bigskip\noindent{{\bf 3. The Tools}}\medskip

The condition $TM=MT$ in (P1) is equivalent to the selection rule
$$M_{\la \mu}\ne 0\Rightarrow (\la+\rho)^2\equiv(\mu+\rho)^2\
\quad({\rm mod}\ 2n)\ .\eqno(3.1)$$
The other (P1) condition $SM=MS$   is much more subtle, as we will see.

The matrix $S$ obeys an important symmetry. Its entries
$S_{\la\mu}$ lie in some cyclotomic extension $\Q(\zeta_N)$ of $\Q$, where
$\zeta_N=\exp[2\pi \i/N]$, so consider any
Galois automorphism $\si\in {\rm Gal}(\Q(\zeta_N)/\Q)\cong\Z_N^\times$. It will
obey $\si\zeta_N=\zeta_N^\ell$ for some $\ell$ coprime to $N$. Then
$$\si(S_{\la\mu})=\epsilon_\si(\la)\,S_{\la^\si\mu}\eqno(3.2a)$$
for some permutation $\la\mapsto \la^\si$ of $P_+$, and some signs 
$\epsilon_\si:P_+\rightarrow\{\pm 1\}$. 
Equation (3.2a) immediately implies [5,15]
$$M_{\la \mu}=\epsilon_\si(\la)\,\epsilon_\si(\mu)\,M_{\la^\si,\mu^\si}
\eqno(3.2b)$$
valid for any $\si$, any $\la,\mu\in P_+$, and any physical invariant $M$.
Positivity (P2) then implies, for all $\sigma$, the selection rule
$$M_{\la \mu}\ne 0\quad\Longrightarrow\quad\epsilon_\sigma(\la)=
\epsilon_\sigma(\mu)\ .\eqno(3.2c)$$
The selection rules (3.1) and (3.2c) are the most important ingredients 
(though there are others) of
step (1) described in section 1. We will use (3.2c) in section 5.

For a positive invariant $M$, define 
$$\eqalignno{\J_L(M)=&\,\{J\in \s_{sc}\,|\,M_{J0,0}\ne 0\}\ ;&(3.3a)\cr
\p_L(M)=&\,\{\la\in P_+\,|\,\exists \mu\in P_+\ {\rm such\ that}\ 
M_{\la \mu} \ne 0\}\ ;&(3.3b)\cr}$$
and define $\J_R(M)$ and $\p_R(M)$ similarly (using the other subscript of 
$M$). Call $\la\in P_+$ a {\it fixed-point} of $\J\subset \s_{sc}$ if 
the cardinalities satisfy $\|\J\la\|<\|\J\|$. 
For example, $\la$ is a fixed-point of $\J_b$ iff $\la_0=\la_1$.
For any $\J\subset \s_{sc}$, define
$$\p(\J)\eqde\{\la\in P_+\,|\,Q_J(\la)\in\Z\ \quad
\forall J\in\J\}\ .\eqno(3.3c)$$

The following elementary lemma tells us
how $\J_L(M)$ and $\J_R(M)$ affect the other entries of $M$.

\medskip\noindent{{\smcap Lemma 3.1}}.\quad[9] (a) {\it Let $M$ be any physical 
invariant, 
and $J,J'\in \s_{sc}$. Then the following statements are equivalent:}

\qquad\qquad\qquad{(i)} $M_{J0,J'0}\ne 0$;

\qquad\qquad\qquad{(ii)} $M_{J0,J'0}=1$;

\qquad\qquad\qquad{(iii)} {\it for any $\la,\mu\in P_+$, if 
$M_{\la \mu}\ne 0$ then $Q_J(\la)\equiv Q_{J'}(\mu)$ (mod $1$)};

\qquad\qquad\qquad{(iv)} $M_{J\la,J'\mu}=M_{\la \mu}$ {\it for all} $\la,\mu
\in P_+$.

\item{{(b)}} {\it Let $M$ be any positive invariant satisfying 
$$M_{0 \mu}=\sum_{J\in\J_R}\delta_{\mu,J0},\qquad M_{\la 0}=
\sum_{J\in\J_L}\delta_{\la,J0}\ ,\eqno(3.4)$$
for some $\J_L,\J_R\subseteq\s_{sc}$. Then} 


\qquad\qquad\qquad(i) $\J_L$ {\it and $\J_R$ are groups and} $\|\J_L\|
=\|\J_R\|$.

\qquad\qquad\qquad(ii) $\p_L(M)=\p(\J_L)$ {\it and} $\p_R(M)=\p(\J_R)$.

\item{(c)} {\it  Let $M$ be any physical invariant satisfying $(3.4)$, for $\J_L
=\J_R=\J$. Suppose no $\la\in \p(\J)$ is a fixed-point of $\J$. Then
there is a permutation $\pi$ of the $\J$-orbits $\p(\J)/\J$ such that
$$\eqalignno{M_{\la \mu}=&\,\left\{\matrix{1&{\rm if}\ 
\mu\in\pi(\J\la)\cr 0&{\rm otherwise}\cr}\right.\ ,&(3.5a)\cr
S_{\la\mu}=&\,S_{\pi\la,\pi\mu}\ ,&(3.5b)\cr}$$
valid for all $\la,\mu\in\p(\J)$ (the other entries of $M$ all vanish).}\medskip

For example, any \ade\ $M$ will obey
(3.4) with $\J_L=\J_L(M)$ and $\J_R=\J_R(M)$.
Note that a special case of Lemma 3.1(c) is that automorphism invariants
(2.1a)  are permutation matrices (2.1b). Lemma 3.1
is a special case of Lemmas 3.1 and 3.2(b) in [9]. In particular,
an analogue of 3.1(c) holds even if $\J_L\ne\J_R$ and if there are fixed-points,
but this simpler case is sufficient for our purposes.

\bigskip\noindent{\bf 4. The level 2 physical invariant
classifications}\medskip

\noindent{{\it 4.1. The $B_{r,2}$  physical invariant
classification.}} The orthogonal algebras at level 2 
behave essentially as if they were rank 1, so the ``brute-force'' approach
of Cappelli-Itzykson-Zuber [2] can be modified to yield an efficient attack
on the problem. In particular, ``unfolding'' $\p_b$ and $\p_v$ puts us into
the familiar terrain of $U_1$ at level $n$($=2r+1$ or $2r$, respectively).
We can explicitly find an integral basis for its commutant, and this
basis quickly solves our classification problem. Three differences between 
our approach and that of [2] are that: (i) for $B_{r,2}$, our commutant 
is for the subgroup $\Gamma_\theta=\langle S,T^2\rangle$ and not the full
modular group $\Gamma=\langle S,T
\rangle$; (ii) folding here preserves positivity; (iii) our approach
applies directly only to the subsets $\p_b$ and $\p_v$ of $P_+(B_{r,2})$ and
$P_+(D_{r,2})$.

We consider first $B_{r,2}$. Recall the notation introduced at the beginning
of section 2.2.
All $S$- and $T$-matrix entries for $B_{r,2}$ are obtained from
$$\eqalignno{S_{00}=&\, {1\over 2}S_{0\ga^a}={1\over\sqrt{n}}S_{0,
J^j\L_r}={1\over 2\sqrt{n}}\ , &(4.1a)\cr
S_{\L_r\L_r}=&\,S_{J\L_r,J\L_r}=-S_{\L_r,J\L_r}={1\over 2}\ ,&(4.1b)\cr
S_{\ga^a\ga^b}=&\,{2\over \sqrt{n}}\cos{2\pi ab\over n}\ ,&(4.1c)\cr
S_{\L_r\ga^a}=&\,S_{J\L_r,\ga^a}=0\ ,&(4.1d)\cr
(\ga^a+\rho)^2=&\,\rho^2+na-a^2\ ,&(4.1e)\cr
(J^j\L_r+\rho)^2=&\,\rho^2+jn+{r\over 4}n\ ,&(4.1f)}$$
for each $a,b\in\{1,\ldots,r\}$ and $j\in\{0,1\}$. The missing values $S_{J0,*}$
can be obtained from (4.1a) by (2.2a), and all other $S$ entries come from $S$
being symmetric. These expressions (4.1a) - (4.1d) immediately follow from 
the calculations leading to rank-level duality, 
and they can also be found in [13]. Note from (4.1e) that 
$M_{\ga^a \ga^b}\ne 0$ requires $a^2\equiv b^2$ (mod $n$). A curiousity 
of this $S$-matrix is that it is essentially the character table for
the dihedral group $D_n$.

The remainder of this subsection is devoted to the proof of Thm.\ 2.1.
We will accomplish this by ``unfolding''. In particular,
define $\widetilde{\p}_n=\Z/n\Z$, and
$$\eqalignno{\tilde{S}_{ab}=&\,{1\over\sqrt{n}}\exp[2\pi \i ab/n]&(4.2a)\cr
\tilde{T}^2_{ab}=&\,\delta_{a,b}\,\exp[2\pi\i a^2/n]&(4.2b)}$$
for all $a,b\in\widetilde{\p}_n$. Then it is easy to see directly that there
is a bijection between all physical invariants $M$ of $B_{r,2}$ with
$\J_L(M)=\J_R(M)=\J_b$, and all nonnegative integral matrices $\widetilde{M}$
obeying:\smallskip

\item\item{($\widetilde{P1}$)} $\widetilde{M}$ commutes with $\tilde{S}$ and 
$\tilde{T}^2$;

\item\item{($\widetilde{P2}$)} $\widetilde{M}_{\pm a,b}=\widetilde{M}_{a,\pm b}=
\widetilde{M}_{a b}$, for any $a,b\in\tilde{\p}_n$;

\item\item{$(\widetilde{P3})$} $\widetilde{M}_{00}=4$, and
$\widetilde{M}_{ab}\in 2\Z$ if either $a=0$ or $b=0$.\smallskip

\noindent{Precisely,} the bijection is given by
$$M_{\ga^a\ga^b}=\widetilde{M}_{ab}\,\left\{\matrix{1&{\rm if\ both}\ a\ne 0\
{\rm and}\ b\ne 0\cr {1\over 4}&{\rm if}\ a=b=0\cr {1\over 2}&
{\rm otherwise}}\right. \eqno(4.2c)$$

\medskip\noindent{\smcap Lemma 4.1}.\quad (a) {\it A basis for the vector space 
$\tilde{\cal V}_n$ of all matrices commuting with $\tilde{S}$ and 
$\tilde{T}^2$, is provided by the set of matrices $\tilde{\B}(d,\ell)$:
$$\tilde{\B}(d,\ell)_{a,b}=\left\{\matrix{1&{\rm if}\ n|da,\ {\rm and}\
b\equiv a\ell\ ({\rm mod}\ d)\cr 0&{\rm otherwise}}\right. \eqno(4.3)$$
{}\qquad where $d|n$, $n|d^2$, and $1\le \ell\le {d^2\over n}$ obeys $\ell^2
\equiv 1$ (mod ${d^2\over n}$).}

\item{(b)} {\it Any integral positive invariant $M$ of $B_{r,2}$ with 
$M_{J0,0}=M_{0,J0}=M_{00}$ can be written as a sum of 
various $\B(d_1,\ell_1|d_2,\ell_2)$.}

\medskip There is a natural geometric interpretation of the $\tilde{\B}(d,
\ell)$  in terms of self-dual lattices, and indeed that interpretation
is the most convenient description of the physical invariant classification
for $U_1\oplus\cdots \oplus U_1$ (see [8]).

\medskip\noindent{\it Proof of Lemma 4.1.}\quad
It is straightforward to verify that the matrices $\tilde{\B}(d,\ell)$ commute
with $\tilde{S}$ and $\tilde{T}^2$. Also, they are all distinct and can be
counted, and we find their number equals the number of divisors of $n$.

These $\tilde{\B}(d,\ell)$ possess an important property:
given any matrix $\tilde{\B}(d,\ell)$, we can find an
index $(i,j)$ such that $\tilde{\B}(d',\ell')_{i,j}\ne 0$ iff $\tilde{\B}(d',
\ell')=\tilde{\B}(d,\ell)$. To see this, choose $\ell_0\equiv \ell$ (mod $d^2/n$)
so that $1\le \ell_0\le n$ and $\ell_0^2\equiv 1$ (mod $n$) -- there may be
more than one $\ell_0$ corresponding to a given $\ell$. Then $\B(d,\ell)=
\B(d,\ell_0)$. The reason $\ell_0$ is more convenient than $\ell$ is that
the $\ell_0^2$ condition means (i) $\ell_0\equiv \pm 1$ (mod $p^a$) whenever
$p^a$ divides $n$, and (ii) together these signs uniquely determine $\ell_0$.
Let $m$ be any integer for which both $m$ and $2\ell_0+m{d^2\over n}$ are coprime
to $n$ ($m$ exists, by the Chinese Remainder Theorem). Then we have
$$\tilde{\B}(d',\ell')_{{n\over d},\ell_0{n\over d}+md}=\left\{\matrix{1&{\rm
if}\ d=d'\ {\rm and}\ \ell_0\equiv\ell'\ ({\rm mod}\ {d^2\over n})\cr
0&{\rm otherwise}}\right.\ ,\eqno(4.4)$$
for any $\tilde{\B}(d',\ell')$ (Proof: The top equation is clear from (4.3);
to see the bottom equation, suppose $\tilde{\B}_{{n\over d},\ell_0{n\over d}
+md}=1$. Then (4.3) says
$d$ must divide $d'$. If $p^a$ divides $dd'/n$ for $a>0$, then $p^a|n$ so
$p$ will be coprime to $m$ and $2\ell_0+md^2/n$, and hence by (4.3) and (i),
$\ell_0\equiv\ell_0'$ (mod $p^a$) and $p^a|{d^2\over n}$. This means by (ii)
that $\ell_0\equiv \ell_0'$ (mod ${dd'\over n}$) and ${dd'\over n}$ divides
${d^2\over n}$. Hence $d=d'$ and $\ell=\ell'$.).

An immediate consequence of (4.4) is that the $\tilde{\B}(d,\ell)$ are
linearly independent. Thus to conclude the proof of part (a), it suffices
to show that dim $\tilde{\cal V}$ is at most the number of divisors of $n$.

We will now follow the proof of Thm.\ 2 
in [5] (which in turn is based on the argument of [2]).
For each $u,u'\in\widetilde{\p}_{2n}=\Z/2n\Z$,
define an $n\times n$  matrix $\{u,u'\}$ by
$$\{u,u'\}_{a,b}=\delta^{(n)}_{a,u+b}\exp[\pi\i\,(u+2b)\,u'/n]\eqno(4.5a)$$
for all $a,b\in\widetilde{\p}_n$, where $\delta^{(n)}_{x,y}$ equals zero unless
$n$ divides $x-y$, when it equals 1. SL$_2(\Z)$ acts on the right by
$\{u,u'\}\left(\matrix{a&b\cr c&d\cr}\right)=\{au+cu',bu+du'\}$. Note that
$$\eqalignno{\tilde{T}^2\{u,u'\}\tilde{T}^2{}^*=&\,\{u,u'\}\left(
\matrix{1&2\cr 0&1}\right)&(4.5b)\cr
\tilde{S}\{u,u'\}\tilde{S}^*=&\,\{u,u'\}\left(\matrix{0&1\cr-1&0}
\right)&(4.5c)\cr}$$
The index of the subgroup $\Gamma_\theta=\langle \left(\matrix{0&1\cr -1&0}
\right),\,\left(\matrix{1&2\cr 0&1}\right)\rangle$ in the modular group
$\Gamma={\rm SL}_2(\Z)$ is 3, with left cosets $\Gamma_\theta$, 
$\left(\matrix{1&1\cr
0&1}\right)\Gamma_\theta$, and $\left(\matrix{1&0\cr 1&1}\right)\Gamma_\theta$.
Certainly $\tilde{\cal V}$ is spanned by the orbit sums 
$N_{\theta}(u,u')\eqde\sum_{g\in\Gamma_\theta}\{u,u'\}g$, but it is more
convenient
to work over the orbit sums $N(u,u')\eqde\sum_{g\in\Gamma}\{u,u'\}g
=(1+(-1)^{u}+(-1)^{u'})\,N_\theta(u,u')$. To see this relation, note that
$\left(\matrix{1&1\cr 0&1}\right)\equiv\left(\matrix{1&2\cr 0&1}\right)^{
{n+1\over 2}}$ and
$$\left(\matrix{1&0\cr1&1}\right)\equiv (\left(\matrix{0&1\cr-1&0}\right)\,
\left(\matrix{1&2\cr 0&1}\right)\,\left(\matrix{0&1\cr-1&0}\right)^{-1})^{
{n-1\over 2}}\ ,$$
both taken (mod $n$). Now for any $u_1,u_2\in\widetilde{\p}_{2n}
$, $N(u_1,u_2)=N(d,0)$ where $d={\rm gcd}(n,u_1,u_2)$.

Thus the number of linearly independent $N(u_1,u_2)$ (and hence the
dimension of $\tilde{\cal V}$) is at most the number of divisors of $n$, and
we are done part (a).

Part (b) is an immediate consequence of Lemma 3.1(b), part (a) and (4.4):
unfolding $M$ into $\widetilde{M}$ via (4.2c), we can write $\widetilde{M}$ as a sum of
$\widetilde{M}_{00}=
4M_{00}$ matrices $\tilde{\B}(d_i,\ell_i)$'s; by $(\widetilde{P2})$,
these can be paired up so that $d_i=d_j$ and $\ell_i=-\ell_j$;
finally, arbitrarily pair up these pairs (possible, since 4 divides $4M_{00}$)
and refold, and we obtain $\B(d_i,\ell_i|d_i',\ell_i')$.\quad\QED\medskip

Lemma 4.1 immediately implies both the $(\widetilde{P1})$ - $(\widetilde{P3})$
classification and, more importantly, the classification of all physical
invariants $M$ of $B_{r,2}$ with $\J_{L,R}(M)=\J_b$. These are precisely
the $\B(d_1,\ell_1|d_2,\ell_2)$ collected in Thm.\ 2.1(b).

Next, consider the possibility that $\J_L(M)=\J_R(M)=\{id.\}$,    
but $M_{\la_s\la_n}=M_{\la_n\la_s}=0$ for
any {\it spinor} $\la_s\in\J_b\L_r$ and any {\it nonspinor} $\la_n\in\p_b$.
By Lemma 3.1(a), this means $M_{J0,J0}=1$, and also
$M_{J\L_r,J\L_r}=M_{\L_r\L_r}$, $M_{0\ga^a}=M_{J0,
\ga^a}$, etc. Of course (2.2b) and (3.1) tell us $M_{J\L_r,\L_r}=M_{\L_r,J\L_r}
=0$. Computing
$SM=MS$ at $(\L_r,0)$ gives us $M_{\L_r\L_r}
=1$, fixing all entries of $M$ involving spinors. Finally,
note that $M'=M\,\I[J]$ is also a physical invariant of $B_{r,2}$ but
with $\J_L(M')=\J_R(M')=\J_b$. Hence $M'=\B(d_1,\ell_1|d_2,\ell_2)$ for some
$d_i,\ell_i$. Since every $M'_{\ga^a\ga^b}$ will be even (for $a,b>0$), 
we must have $d_1=d_2$ and $\ell_1=\ell_2$. This 
fixes all entries of $M$, and we get $M=\B(d_1,\ell_1)$.

 From the calculations (4.1a) and the Galois selection rule
 (3.2c), this concludes the proof whenever $n$ is not a perfect square.
 So consider now $\sqrt{n}\in\Z$. We may assume $M_{J0,J0}=0$
 and without loss of generality that $M_{J0,0}=0$.
Then $T$-invariance (3.1) and (4.1e),(4.1f) say that for any $\nu\ne \mu^r$,
$$M_{\mu^r\nu}=M_{\nu\mu^r}=0\ .\eqno(4.6a)$$
$MS=SM$ evaluated at $(0,\mu^r)$, $(\mu^r,0)$ and $(J0,\mu^r)$,
$(\mu^r, \ga^a)$, and $(\mu^r,\la^r)$ gives us
$$\eqalignno{M_{\mu^r\mu^r}=&\,1-M_{0,J0}-M_{0\la^r}=1-M_{\la^r
0}=M_{J0,\la^r}&(4.6b)\cr
M_{\la^r\ga^a}=&\,M_{0\ga^a}-M_{J0,\ga^a}&(4.6c)\cr
M_{\la^r\la^r}=&\,M_{0\la^r}&(4.6d)\cr}$$
for all $a>0$. Note that (3.1) and (4.1e),(4.1f) say $M_{\la^r\ga^a}=
M_{\ga^a\la^r}=0$ unless $\sqrt{n}$ divides $a$. So comparing
$MS=SM$ at $(\la^r,0)$ and $(\la^r,\ga)$ for any $\ga\in\c$,
and using (4.6c), gives
$$M_{\la^r\ga}=1+M_{\la^r0}-M_{\la^r\la^r}\ .\eqno(4.6e)$$
By (4.6b), there are three possibilities:

\smallskip
\item\item{(i)} $M_{0,J0}=M_{\la^r0}=1$ and $M_{\mu^r\mu^r}=M_{*,\la^r}=0$\ ;

\item\item{(ii)} $M_{0,J0}=M_{\mu^r\mu^r}=M_{J0,\la^r}=0$ and
$M_{0\la^r}=M_{\la^r0}=1$\ ;

\itemitem{(iii)} $M_{0,J0}=M_{0\la^r}=M_{\la^r0}=0$ and $M_{\mu^r
\mu^r}=M_{J0,\la^r}=1$\ .\smallskip

In possibility (i), consider the product $\I[J]\, M$: it will be a physical
invariant and so by Lemma 4.1(b) will equal some $\B(d,\ell|d,\ell)$.
Hence for each $\ga^a$, 
$$M_{0\ga^a}+M_{J0,\ga^a}=(\I[J]\, M)_{0\ga^a}=(\I[J]\, M)_{
\ga^a 0}=2\,M_{\ga^a 0}$$
will equal either 0 or 2. Together with (4.6c) and (4.6e), we get $M_{\la^r
\ga}=M_{0\ga}=2$, $M_{\ga 0}=1$ and $M_{J0,\ga}=0$ for $\ga\in\c$.
Hence we must have $\I[J]\, M=\B(\sqrt{n},1|\sqrt{n},1)$, and we can read
off the remaining entries: $M_{\ga\ga'}=2$ or 0, depending on whether
or not both $\ga$ and $\ga'$ lie in $\c$. We thus obtain $M=\B^{iii}$.

For possibility (ii), use $M\,\I[J]=\B^{iii}$ and $\I[J]\, M=(\B^{iii})^T$, in
order to show $M=\B^{ii}$. Similarly, for possibility (iii) we find $M=
\B^{i}$.

\bigskip\noindent{{\it 4.2.\quad The $D_{r,2}$ physical invariant classification.}}
Next we consider $D_{r,2}$. The argument is very
analogous to the $B_{r,2}$ one. Recall the notation introduced
at the beginning of section 2.3. 
 By the usual calculations (e.g.\ writing $(S_{\L_r\L_r}\pm S_{\L_r\L_{r-1}})/
S_{0\L_r}$ as a product of sines/cosines) we get 
$$\eqalignno{S_{00}=&\,{1\over \sqrt{r}}S_{0\L_r}={1\over 2}
S_{0\la^a}={1\over 2\sqrt{n}},&(4.7a)\cr S_{\la^a\la^b}=&\,
{2\over \sqrt{n}}\cos(\pi{ab\over r}),&(4.7b)\cr S_{\la^a\L_r}=&\,
S_{\la^a\L_{r-1}}=0&(4.7c)\cr S_{\L_r\L_r}=&\,S_{\L_{r-1}\L_{r-1}}=
{1\over 4}(1+(-\i)^r),&(4.7d)\cr S_{\L_r\L_{r-1}}=&\,{1\over 4}(1-(-\i)^r),
&(4.7e)
\cr(\la^c+\rho)^2= &\,\rho^2+2rc-c^2,&(4.7f)\cr (J^j\L_r+\rho)^2=&\,
(J^j\L_{r-1}+\rho)^2=\rho^2+jn+{r^2\over 2}-{r\over 4},&(4.7g)\cr}$$
where $a,b\in\{1,2,\ldots,r-1\}$, $c\in\{0,1,\ldots,r\}$ and $j\in\{0,1\}$.
The remaining entries of $S$ are given by (2.2a) and $S=S^T$.
Again we have the curious
relation between this matrix $S$ and the dihedral group $D_{2r}$.

We will ``unfold'' $D_{r,2}$ as we did $B_{r,2}$. Namely,
define $\widetilde{\p}_n=\Z/n\Z$ as before, and write
$$\eqalignno{\tilde{S}_{ab}=&\,{1\over\sqrt{n}}\exp[2\pi \i ab/n]&\cr
\tilde{T}_{ab}=&\,\delta_{a,b}\,\exp[\pi\i a^2/n]&}$$
for all $a,b\in\widetilde{\p}_n$. The difference here is that we are able to
define $\tilde{T}$ rather than merely $\tilde{T}^2$ -- this simplifies
the arguments. It is easy to see directly that there
is a bijection between all physical invariants $M$ of $D_{r,2}$ with both
$J_v\in\J_L(M),J_v\in\J_R(M)$, and all nonnegative integral $n\times
n$ matrices $\widetilde{M}$ obeying:

\item\item{$(\widetilde{P1})$} $\widetilde{M}$ commutes with $\tilde{S}$ and 
$\tilde{T}$;

\item\item{($\widetilde{P2}$)} $\widetilde{M}_{\pm a,b}=\widetilde{M}_{a,\pm b}=
\widetilde{M}_{ab}$, for any $a,b\in\widetilde{\p}_n$;

\item\item{$(\widetilde{P3})$} $\widetilde{M}_{00}=4$, $\widetilde{M}_{ab}\in 4\Z$ for
$a,b\in\{0,r\}$, $\widetilde{M}_{ab}\in 2\Z$ if either $a\in\{0,r\}$ or
$b\in\{0,r\}$.

\noindent{In} fact the bijection is given by
$$M_{\la^a\la^b}=\widetilde{M}_{ab}\,\left\{\matrix{1&{\rm if\ both}\ a,b\not\in
\{0,r\}\cr {1\over 4}&{\rm if\ both}\ a,b\in\{0,r\}\cr 
{1\over 2}&{\rm otherwise}}\right. $$

\medskip\noindent{\smcap Lemma 4.2}.\quad  (a) {\it A basis for the vector space 
$\tilde{\cal V}_n$ of all matrices commuting with $\tilde{S}$ and 
$\tilde{T}$, is provided by the set of all matrices $\tilde{\B}(d,\ell)$
given by $(4.3)$, where here $d|n$, $2n|d^2$, and $1\le \ell\le
{d^2\over n}$ obeys $\ell^2\equiv 1$ (mod ${2d^2\over n}$).}

\item{(b)} {\it Any positive integral invariant $M$ of $D_{r,2}$ with 
$M_{J0,0}=M_{0,J0}=M_{00}$, can be written as a sum of 
various $\D(d_1,\ell_1|d_2,\ell_2)$.}

\medskip Lemma 4.2(a) is a special case of Thm.2 in [5] and follows
from a simplified version of our proof of Lemma 4.1. 
Counting the dimension of $\tilde{{\cal V}}_n$ as in Lemma 4.1(a),
we see that it equals the number of divisors of $n/2=r$.

Lemma 4.2(b) immediately gives us the classification of all $D_{r,2}$
physical invariants with $J_v\in\J_{L,R}(M)$: they are the $\D(d_1,\ell_1|d_2,
\ell_2)$ given in Thm.\ 2.2(b).

Next, suppose both $J\not\in\J_{L}(M)$ and $J\not\in\J_R(M)$, but $M_{\la_s
\la_n}=M_{\la_n\la_s}=0$ for any spinor $\la_s\in\{C_1^iJ^j\L_r\}$, and
any non-spinor $\la_n\in\p_v$. By Lemma 3.1(a), this means $M_{J0,
J0}=1$. Comparing $SM=MS$ at $(\la_s,J'0)$ for each choice of $\la_s\in\{
\L_{r-1},\L_r\}$ and $J'\in\J_s$, we find that either: 

\item{--} (replacing $M$ if necessary by its conjugation $M\, C_1$) 
$M_{J_s0,
J_s0}=M_{\L_r\L_r}=M_{\L_{r-1}\L_{r-1}}=1$, $M_{J_s0,J_c0}=
M_{\L_r\L_{r-1}}=M_{\L_{r-1}\L_r}=0$, and $\J_{L,R}(M)=\{id.\}$;

\item{--} (replacing $M$ with some  $C_1^i\, M\, C_1^j$)
$M_{\L_r\L_r}=2$, $\L_{r-1}\not\in\p_{L,R}(M)$, and $\J_{L,R}(M)=\J_s$; or

\item{--} (replacing $M$ with some $C_1^i\,M\,C_1^j$ and if necessary
transposing) $\J_L(M)=\J_s$ and $\J_R(M)=\{id.\}$.

\noindent The third possibility is eliminated by evaluating $SM=MS$ at (0,0):
the left side is an even multiple of $S_{00}$, while the right side is an
odd multiple of $S_{00}$.
The remaining entries of $M$ in the first two cases are fixed by Lemma
4.2(b) and the projection $M\mapsto M\, \I[J]$. We find that in either
case, $M$ (appropriately conjugated) equals one of the $\D(d,\ell)$
of Thm.2.2(a).

This concludes the classification of all physical invariants whenever $r$
is not a perfect square, by the Galois argument (3.2c) applied to (4.7a),
or whenever 4 does not divide $r$, by T-invariance (3.1). So consider now $r$
a perfect square, $4|r$, $M_{J0,J0}=0$ and without loss of generality
$M_{J0,0}=0$. $T$-invariance also says that
spinors cannot couple to $\J_d0$. Recall the definition of ${\cal C}_j$
given in subsection 2.3. Then $T$-invariance says that
$\la^a$ can couple to $\J_d0$  only for
$\la^a\in{\cal C}_0\cup {\cal C}_4$, $\la^b$
can couple to $\L_{r-1}$ or $\L_r$ only for $\la^b\in{\cal C}_1$, and $\la^c$
can couple to $J\L_{r-1}$ or $J\L_r$ only for $\la^c\in{\cal C}_3$.

Lemma 4.2(b) tells us that  either $\I[J]\, M$ (if
$M_{0,J0}=1$) or $\I[J]\, M\,\I[J]$ (if $M_{0,J0}=0$)
equals $\D(d_1,\ell_1|d_2,\ell_2)$ for $d_i,\ell_i$ as in Thm.\ 2.2(b).
Because the $(2d_i,0)$ and $({r\over d_i},\ell_i{r\over d_i})$ entries of
the product will be even, we find that $d_1=d_2\eqde d$ and $\ell_1=\ell_2\eqde
\ell$.

For any choice of $\la_s\in\{\L_{r-1},\L_r\}$, 
$MS=SM$ evaluated at $(\la_s,\la)$  gives
$${2\over \sqrt{r}}\sum_{\nu\in{\cal C}_1}M_{\la_s\nu}=
M_{\la_s\la'}={\pm 1\over 2}(M_{0\la_{\pm}}-M_{J0,\la_{\pm}}
+s(\la_s)M_{J_s0,\la_{\pm}}-s(\la_s)M_{J_c0,\la_{\pm}})$$
where $\la_+\in{\cal C}_0$, $\la_-\in{\cal C}_4$, $\la'\in{\cal C}_1$, and
$s(\la_s)=\exp[2\pi\i\,Q_s(\la_s)]\in\{\pm 1\}$.
The first inequality says that, for fixed $\la_s\in \{\L_{r-1},\L_r\}$,
$M_{\la_s\la'}$
is independent of $\la'\in{\cal C}_1$ (call this value $\M_L(\la_s)$). Evaluating
$SM=MS$ at $(\la_s,\mu)$ for $\mu\in{\cal C}_3$ shows $M_{J\la_s,\mu}$ is also
constant and equals $\M_L(\la_s)$.

Comparing $\la_s\in\{\L_{r-1},\L_r\}$ in the second equality, we get
$M_{0\la}>M_{J0,\la}$ for $\la\in{\cal C}_0$, and $M_{0\la}<
M_{J0,\la}$ for $\la\in{\cal C}_4$ (equalities here would mean no spinors
lie in $\p_L$, hence would contradict $J\not\in\J_L(M)$). Choosing
$\la=\la^{2\sqrt{r}}$ here then forces
$d=\sqrt{r}$, hence $\ell=1$. Since $J_s\in\J_{L,R}(\D(\sqrt{r},1|\sqrt{r},1))$,
we know $J_s\in\J_{L,R}(M)$, provided we conjugate $M$ appropriately. One consequence of this is that 16 must
divide $r$, since if $4\|r$, then ${\cal C}_1\cap\p_s=\emptyset$. Another
consequence is that $\L_{r-1}\not\in\p_{L,R}(M)$, and so ${\cal M}_L(\L_r)
$ must be positive.

Now consider those $M$ with $M_{0,J0}=1$.  For $\la\in{\cal C}_0$, we know
${\cal M}_L(\L_r)=M_{0\la}-M_{J0,\la}$
must be positive and independent of $\la$. Since $M_{0
\la}+M_{J0,\la}= 2$, the only possibility is that
$M_{0\la}={\cal M}_L(\L_r)=2$ and $M_{J0,\la}=0$.
Similarly, for $\la\in{\cal C}_4$ we find $M_{J0,\la}=2$ and $M_{0\la}
=0$. This determines all entries of $M$, and we obtain ${\cal D}^{ii}$.

When $M_{0,J0}=0$ the equality $M\,\I[J]=\D^{ii}$ fixes most entries of $M$.
$M_{J^i\L_r,J^i\L_r}=1$
is forced by evaluating $MS=SM$ at $(J^i\L_r,0)$. 
We thus obtain ${\cal D}^{i}$.

\bigskip\bigskip\noindent{{\bf 5. The level 3 physical invariant
classification}}\bigskip

Write $0(M)$ for the set of all weights coupled to $0$ (i.e.\ all weights
$\la$ obeying either $M_{0\la}\ne 0$ or $M_{\la 0}\ne 0$). We want
to show that $0(M)\subset\s_{sc}0$ for any physical invariant $M$ of
$B_{r,3}$ or $D_{r,3}$ (except for the $D_{7,3}$ exceptionals).

Consider first $B_{r,3}$. Write $n=2r+2$, $\ga^i\eqde \L_i$ for $i<r$, $\ga^r
\eqde2\L_r$, $\mu^i\eqde \L_i+\L_r$ for $1\le i<r$, 
and $\mu^r\eqde 3\L_r$. There are $3r+4$ weights in $P_+$: 0, $3\L_1=J_b0$,
$\L_r$, $2\L_1+\L_r=J_b\L_r$, $\ga^i$, $J_b\ga^i$, and $\mu^i$. The norms are
$$\eqalignno{(J_b^i\L_r+\rho)^2=&\,
\rho^2+{r^2\over 2}+{r\over 4}+2ni&(5.1a)\cr
(J_b^i\ga^a+\rho)^2=&\,\rho^2+a\,(2r+1-a)+ni&(5.1b)\cr
(\mu^a+\rho)^2=&\,\rho^2+{r\over 4}(2r+1)+a\,(n-a)\ ,&(5.1c)\cr}$$
for $i\in\{0,1\}$ and $1\le a\le r$, while the q-dimensions are given by
$$\eqalignno{\D(\ga^a)=&\,{\sin(\pi\,(2a+1)/2n)\over \sin(\pi/2n)}&(5.2a)\cr
\D(\mu^a)=&\,\sqrt{2}\,{\sin(\pi\,(r+1-a)/n)\over \sin(\pi/2n)}\ .&(5.2b)\cr}$$
The main use of q-dimensions here is for reading off the Galois
parities: $\eps_\si(\la)\,\eps_\si(0)={\rm sgn}[\si\D(\la)]$ (see (3.2)).
The automorphisms $\si_\ell$ are parametrised by integers
$0<\ell<4n$ coprime to $n$. For example, $\si_\ell\sqrt{2}=\sqrt{2}$
if $\ell\equiv\pm 1$ (mod 8), otherwise it equals $-\sqrt{2}$.

The first step is to show that no spinor $\la$ (i.e.\ $\la\in P_+$ with
$\la_r$ odd) can couple to $0$. Equation (5.1a) tells us $\L_r,
J_b\L_r\not\in 0(M)$. Suppose for contradiction that $\mu^a\in 0
(M)$ for $a$ odd. Equation (5.1c) says 4 divides $r$. Take $\ell=n-1$: then
(5.2b) says $\eps_\si(\mu^a)\,\eps_\si(0)=-1$, which contradicts the
Galois selection rule (3.2c). Next, suppose instead $\mu^a\in 0(M)$
for $a$ even. Equation (5.1c) says 8 divides $r$. Then checking each of the
4 cases (${a\over 2}$ even/odd; $a\le{r\over 2}$ or $a>{r\over 2}$), we find
that one of $\ell={n\over 2}\pm 2$ will violate the Galois selection rule.

Thus the only possible weights coupling to $0$ are in $\J_b\ga^a$ and
$\J_b0$. The Galois 
selection rule (3.2c) for them reduces precisely to that of $A_{1,4r+2}$,
which was solved in Lemma 5 of [6]. Incidently, the proof of that lemma
 could have been simplified {\it enormously} by rewriting the parity
trigonometrically: 
$$\epsilon'_\ell(a\L_1')\,\epsilon'_\ell(0')={\rm sgn}(\cos(\pi\ell a/m)
-\cos(\pi\ell\,(a+2)/m))$$
where primes denote the quantities in $A_{1,m-2}$. Anyways,
what we find is that for $r\ge 3$,
no $\ga^a$ can satisfy the Galois selection rule (there is however
a solution
for $r=2$, corresponding to the $C_{2,3}$ exceptional). Thus every physical
invariant for $B_{r,3}$ is an \ade.

Moreover, $(J_b0+\rho)^2\equiv \rho^2+n$ (mod 2$n$), and hence any \ade\
will be an automorphism invariant. These were classified in [11], and
we find that the only ones for $B_{r,3}$ are $M=I$ and
$M=\I[J_b]$.

Incidently, it is particularly easy to classify the $B_{r,3}$ automorphism
invariants: q-dimensions and (3.5b) tell us that $\pi\la\in\J_b\la$, and then
use fusion coefficients or the values $S_{\L_1\la}$ to show that $\pi\L_r
=\L_r$ implies $M=I$, whereas $\pi\L_r=J_b\L_r$ implies $M=\I[J_b]$.

The proof for $D_{r,3}$ is easier. Put $n=2r+1$. Then each of the $4r+8$ 
weights in $P_+$
can be mapped by a simple current $J$ and possibly $C_1$ to one of $\ga^0\eqde
0$,
$\ga^i\eqde\L_i$ for $1\le i\le r-2$, $\ga^{r-1}\eqde\L_{r-1}+\L_r$,
or $\ga^r\eqde 2\L_r$. Norms and q-dimensions are given by
$$\eqalignno{(\ga^a+\rho)^2=&\,\rho^2+a\,(2r-a)&(5.3a)\cr
\D(\ga^a)=&\,{\sin(\pi\,(2a+1)/2n)\over \sin(\pi/2n)}&(5.3b)\cr}$$
with a factor of ${1\over 2}$ on the right side of (5.3b) if $a=r$.
The norms for $C_1^iJ\ga^a$ of course can now be obtained using (2.2b),
while its q-dimension equals that of $\ga^a$. Lemma 5 of [6]
again applies,
and we find that the only possibility for an anomolous coupling with $0$
is $\J_d\ga^5$ at $r=7$. The usual arguments (see e.g.\ [7]) allow us to
construct the exceptional physical invariant, and we obtain either $\E(D_{7,3})$
or $C_1\,\E(D_{7,3})$.

All other $D_{r,3}$ physical invariants will be \ade s. The automorphism
invariants were classified in [11]: we get $M=C_i$ and $C_i\,\I[J_v]$, as well
 as (for $r\equiv 4$ (mod 8)) the additional ones in
$C_i\langle\I[J_s],\I[J_c]\rangle$, for arbitrary
conjugation $C_i$.

The remaining \ade s are found in [10], but we will provide an alternate
 argument
here. Let $M$ be any \ade\ which is not an automorphism invariant.
Hitting $M$ on either side
if necessary with $C_1$, Lemma 3.1(b)(i) and (3.1) tell us that it is sufficient
to consider $\J_L(M)=\J_R(M)=\J_s$, with $8|r$. There are precisely
$r+2$ orbits $\p_s/\J_s$: $\J_s0$, $\J_s(J_v^a\ga^a)$ for $1\le a\le r$, and
$\J_s(2\L_{r-1})$. Note that there are no $J_s$-fixed-points (since the level
is odd) -- this simplifies enormously the argument. This means (Lemma 3.1(c))
that there is a permutation $\pi$ on those $r+2$ orbits such that
(3.5) holds. In particular, the q-dimensions $\D(\la)$ and $\D(\pi\la)$
must be equal. Now (5.3b) tells us $\D(\ga^0)<\D(\ga^1)<\cdots<\D(\ga^{r-1})$
and $\D(\ga^r)=\D(\ga^{{r-1\over 3}})$, 
and hence $\pi(\J_sJ_v^a\ga^a)=\J_s J_v^a\ga^a$ for $a<r$ (to eliminate $\pi(
\J_sJ_v
\ga^{{r-1\over 3}})=\J_s C_1^i\ga^r$, use the fact that $\pi$ must be a symmetry
of the fusion $(J_v\L_1)\stimes\L_{{r-1\over 3}-1}=(J_v\L_{{r-1\over 3}-2})
\splus(J_v\L_{{r-1}\over 3})\splus(\L_{{r-1\over 3}-1})$). Then
$M=\I[J_s]$
if $\pi$ fixes $\J_s(2\L_r)$, and $M=C_1\,\I[J_v]\,\I[J_s]$ if
instead $\pi(\J_s(2\L_r))=\J_s(2\L_{r-1})$.

\bigskip\bigskip\noindent{{\bf 6. Concluding remarks.}} \bigskip

Something very unusual happens for $B_r^{(1)}$ and $D_r^{(1)}$
at level 2, as has been noticed previously in the literature [4,11,16].
Indeed, many of the techniques available for generic algebras
and levels -- most significantly the Galois selection rule (3.2c) -- break down
at $B_{r,2}$ and $D_{r,2}$. This is a symptom of the existence
here of a large family of exceptionals and is our motivation for doing their
 physical invariant classification.
In this paper we also accomplish this classification for
$B_{r,3}$ and $D_{r,3}$ -- they follow quickly from a lemma solving the
Galois selection rule for $A_1^{(1)}$. Thus this paper adds four
 more notches to the surprisingly
barren bedpost representing families of $X_{r,k}$ for which the physical
invariant classification has been completed. In the process we find infinitely
many new exceptionals for both $B_{r,2}$ and $D_{r,2}$. The only level 3
exceptionals occur at $D_{7,3}$ and $B_{2,3}$ (which more properly should
be written $C_{2,3}$).

Our explanation for the rich structure of level 2 physical invariants should
be clear from the argument of section 4: rank-level duality relates
$B_{r,2}$ and $D_{r,2}$ to $U_{1,2r+1}$ and $U_{1,2r}$, respectively,
and $U_{1,n}$ has a known rich family of physical invariants [8]. An alternate
explanation is offered in [16], using the $c=1$ orbifolds $SO(N)_1\times
SO(N)_1/SO(N)_2$.

The only other low-level classifications, for any of the algebras, which are 
important are $B_{r,k}$ and $D_{r,k}$ for $k=4$ and 8, and to a lesser extent 
all other $k\le 6$ for these orthogonal algebras and the trivial case $A_{r,1}$.
The  reason again is rank-level duality: it breaks down or at least takes a
different form for these algebras and levels.  The 
reason the physical invariant classification for $B_{r,8}$ and $D_{r,8}$
would be interesting is that rank-level duality associates to it
the very special algebra $D_4^{(1)}$, and 
$D_4$ triality is already known to give families of $B_{r,8}$ and $D_{r,8}$
exceptionals [18]. 
$B_{r,4}$ and $D_{r,4}$ are interesting because their Galois selection rules
can be solved but have many solutions -- this is usually a sign of 
exceptional chiral extensions (hence exceptional physical invariants).
It should be possible to do these classifications with our current understanding.

The $C_{2,k}$ classification should be straightforward, and would also
imply the $C_{r,2}$, $B_{r,5}$ and $D_{r,5}$ classifications.
All $C_{2,k}$ physical invariants are known for $k\le 500$, and exceptionals
appear only at $k=3,7,8,12$.
Much more valuable, but more difficult, would be the $G_{2,k}$ classification.
Its only exceptionals for $k\le 500$ appear at $k=3,4$. A very safe conjecture
is that there are no new $C_{2,k}$ and $G_{2,k}$ exceptionals.

\bigskip\noindent{{\it Acknowledgments.}}\quad
This research was supported in part by NSERC. Part of this paper was written
at Feza G\"ursey  Institute in Istanbul, and I thank it
for its warm hospitality.\bigskip

\bigskip \noindent{{\bf References}} \bigskip

\item{1.} Bernard, D.: String characters from Kac-Moody 
automorphisms.  Nucl.\ Phys.\ {\bf B288}, 628-648 (1987)

\item{2.} Cappelli, A., Itzykson, C., Zuber, J.-B.: The 
A-D-E classification of $A_1^{(1)}$ and minimal conformal field theories.
 Commun.\ Math.\ Phys.\ {\bf 113}, 1-26 (1987)

\item{3.} Evans, D.E., Kashahigashi, Y.: Quantum Symmetries on
Operator Algebras. Oxford: Oxford University Press, 1998

\item{4.} Fuchs, J., Schellekens, A.N., Schweigert, C.:  
Galois modular invariants of WZW models. 
Nucl.\ Phys.\ {\bf B437}, 667-694 (1995)

\item{5.} Gannon, T.: WZW commutants, lattices, and level-one
partition functions. 
Nucl.\ Phys.\ {\bf B396}, 708-736 (1993)

\item{6.} Gannon, T.: The classification of SU(3) modular
invariants revisited.
 Annales de l'I.H.P.: Phys.\ Th\'eor.\ {\bf 65}, 15-56 (1996)

\item{7.} Gannon, T.: The level 2 and 3 modular invariant 
partition functions for SU(n). 
 Lett.\ Math.\ Phys.\ {\bf 39}, 289-298 (1997)

\item{8.} Gannon, T.: $U(1)^m$ modular invariants, N=2 minimal models,
and the quantum Hall effect. 
 Nucl.\ Phys.\ {\bf 491}, 659-688 (1997)

\item{9.} Gannon, T.:  Kac-Peterson, Perron-Frobenius, and
the classification of conformal field theories. q-alg/9510026

\item{10.}  Gannon, T.: The \ade s of affine algebras. (in preparation)

\item{11.} Gannon, T., Ruelle, Ph., Walton, M.A.:
Automorphism modular invariants of current algebras. 
 Commun.\ Math.\  Phys.\ {\bf 179}, 121-156 (1996)

\item{12.} Kac, V.G.: Infinite Dimensional Lie algebras, 
3rd edition. Cambridge: Cambridge University Press,  1990

\item{13.} Kac, V.G.,  Wakimoto, M.: Modular and conformal 
constraints in representation theory of affine algebras.
 Adv.\ Math.\ {\bf 70}, 156-236 (1988)

\item{14.} Mlawer, E.J., Naculich, S.G., Riggs, H.A.,
Schnitzer, H.J.:  Group-level duality of WZW fusion coefficients and 
Chern-Simons link observables.  Nucl.\ Phys.\ {\bf B352}, 863-896 (1991)

\item{15.} Ruelle, Ph., Thiran, E.,  Weyers, J.:  
Implications of an arithmetic symmetry
of the commutant for modular invariants. 
 Nucl.\ Phys.\ {\bf B402}, 693-708 (1993)

\item{16.} Schellekens, A.N.: Cloning SO(N) level 2. math.QA/9806162

\item{17.} Schellekens, A.N., Yankielowicz, S.: Modular 
invariants from simple currents. An explicit proof. 
 Phys.\ Lett.\   {\bf B227}, 387-391 (1989)

\item{18.} Verstegen, D.: New exceptional modular invariant
partition functions for simple Kac-Moody algebras. 
 Nucl.\ Phys.\ {\bf B346}, 349-386 (1990)

\end